\documentclass{amsart}
\usepackage{amssymb,latexsym}
\usepackage{euscript}
\usepackage{caption2}
\usepackage[dvips]{color}
\begin{document}
\input epsf.sty

\newtheorem{theorem}{Theorem}[section]
\newtheorem{lemma}[theorem]{Lemma}
\newtheorem{corollary}[theorem]{Corollary}
\newtheorem{proposition}[theorem]{Proposition}

\title[Jordan-H\"older theorem for imprimitivity systems]{Jordan-H\"older theorem for imprimitivity systems and maximal decompositions of rational functions}
\author {M. Muzychuk, F. Pakovich}
\markright{Generalized Jordan-H\"older theorem}

\date{}


\begin{abstract} In this paper we prove several results about the lattice of imprimitivity systems
of a permutation group containing a cyclic subgroup with at most two orbits. As an application we
generalize the first Ritt theorem
about functional decompositions of polynomials, and some other related results.
Besides, we discuss examples of rational functions, related to finite subgroups of ${\rm Aut}(\mathbb C{\mathbb P}^1),$
for which the first Ritt theorem fails to be true.
\end{abstract}

\maketitle
\def\be{\begin{equation}}
\def\ee{\end{equation}}
\def\bs{$\square$ \vskip 0.2cm}
\def\d{{\rm d}}
\def\D{{\rm D}}
\def\I{{\rm I}}
\def\C{{\mathbb C}}
\def\N{{\mathbb N}}
\def\P{{\mathbb P}}
\def\Z{{\mathbb Z}}
\def\R{{\mathbb R}}
\def\ord{{\rm ord}}
\def\f{\EuScript}

\def\e{\eqref}
\def\phi{{\varphi}}
\def\v{{\varepsilon}}
\def\deg{{\rm deg\,}}
\def\Aut{{\rm Aut\,}}
\def\Det{{\rm Det}}
\def\dim{{\rm dim\,}}
\def\Ker{{\rm Ker\,}}
\def\Gal{{\rm Gal\,}}
\def\St{{\rm St}}
\def\exp{{\rm exp\,}}
\def\cos{{\rm cos\,}}
\def\diag{{\rm diag\,}}
\def\GCD{{\rm gcd }}
\def\LCM{{\rm lcm }}
\def\mod{{\rm mod\ }}
\def\c{\circ}
\def\eq{{\rm Eq}}

\def\bp{\begin{proposition}}
\def\ep{\end{proposition}}
\def\bt{\begin{theorem}}
\def\et{\end{theorem}}
\def\be{\begin{equation}}
\def\bee{\begin{equation*}}
\def\la{\label}
\def\l{\lambda}
\def\m{\mu}
\def\ee{\end{equation}}
\def\eee{\end{equation*}}
\def\bl{\begin{lemma}}
\def\el{\end{lemma}}
\def\bc{\begin{corollary}}
\def\ec{\end{corollary}}
\def\br{\noindent{\bf Remark.\ }}
\def\pr{\noindent{\it Proof. }}
\def\note{\noindent{\bf Note. }}
\def\bd{\begin{definition}}
\def\ed{\end{definition}}
\def\qed{$\ \ \Box$ \vskip 0.2cm}
\newcommand{\sg}[1]{\langle{#1}\rangle}
\newcommand{\New}[1]{{\sl #1}}
\newcommand{\bS}{{\bf S}}
\newcommand{\cP}{\mathcal P}
\newcommand{\cE}{\mathcal E}
\newcommand{\lcm}{{\rm lcm}}
\newcommand{\sym}{{\rm Sym}}
\newcommand{\cR}{\mathcal{R}}
\newcommand{\new}[1]{{\emph {#1}}}
\newcommand{\cmnt}[1]{\textcolor{black}{#1}}
\newcommand{\normaleq}{\trianglelefteq}
\newcommand{\LG}{L(G_1,G)}
\newcommand{\core}{{\sf core}}
\newcommand{\change}[1]{\textcolor{red}{#1}}

\section{Introduction} Let $F$ be a rational function with complex coefficients.
The function $F$ is called {\it indecomposable} if
the equality $F=F_1\circ F_2,$
where $F_1\circ F_2$ denotes the superposition $F_1(F_2(z))$
of rational functions $F_1,$ $F_2$,
implies that at least one of the functions $F_1,F_2$ is of degree 1.
A rational function which is not indecomposable is called {\it decomposable}. Any representation $\f F$
of a rational function $F$ in the form \be \la{dec} F=F_1\circ F_2\circ \dots \circ F_r,\ee
where $F_1,F_2,\dots, F_r$ are rational functions,
is called {\it a decomposition} of $F.$
If all $F_1,F_2,\dots, F_r$ are
indecomposable of degree greater than one, then the decomposition $\f F$ is called {\it maximal}.
Two decompositions of a rational function $F$
\be \la{decc} F=U_1\circ U_2\circ \dots \circ U_k \ \ \ \ {\rm and} \ \ \ \
F=V_1\circ V_{2}\circ \dots \circ V_m,\ee maximal or not,
are
called {\it equivalent} if they have the same length (that is $k=m$) and there exist rational functions of degree one
$\mu_i,$ $1\leq i \leq k-1,$ such that
$$U_1=V_1\circ \mu_1, \ \ \
U_i=\mu_{i-1}^{-1}\circ U_i \circ \mu_{i}, \ \ \ 1<i< k, \ \ \ {\rm and} \ \ \ V_k=\mu_{k-1}^{-1}\circ V_k.
$$

In the paper \cite{r1} Ritt described the structure of possible maximal decompositions of
polynomials.
This description can be summarized in the form
of two theorems usually called the first and the second Ritt
theorems (see \cite{r1}, \cite{sch}).
The first Ritt theorem states that for any two maximal decompositions $\f D, \f E$
of a
polynomial $F$
there exists
a chain of maximal decompositions $\f F_i$, $1\leq i \leq s,$ of $F$ such that
$\f F_1=\f D,$ $\f F_s\sim \f E,$ and $\f F_{i+1}$ is obtained from $\f F_i,$
$1\leq i \leq s-1,$ by replacing two successive functions in $\f F_i$ by other functions with the same composition.
This implies in particular that any two maximal decompositions of a polynomial
have the same length.
Below we will call two maximal decompositions $\f D, \f E$ of a rational function $F$
such that there exists a chain as above {\it weakly equivalent}.
This defines an equivalence relation on the set of maximal decompositions of $F$.

The first Ritt theorem reduces the description of maximal decompositions of
polynomials to the description of indecomposable polynomial solutions of the equation
\be \la{-0} A\circ C=B\circ D\ee such that the decompositions $A\circ C$ and $B\circ D$ are non-equivalent,
and the second Ritt theorem states if $A,B,C,D$ is such a solution
then there exist polynomials
$\hat A, \hat B, \hat C, \hat D,$ $\mu_1, \mu_2,$ where $\deg \mu_1=1, \deg \mu_2=1,$ such that
$$A=\mu_1\circ \hat A, \ \  B=\mu_1\circ \hat B,\ \ C=\hat C\circ \mu_2, \ \  D=\hat D \circ \mu_2, \ \ \hat A\circ \hat C=\hat B\circ \hat D,$$ and
up to a possible replacement of
$\hat A$ by $\hat B$ and $\hat C$ by $\hat D$ either
$$\hat A\circ \hat C\sim z^n \circ z^rR(z^n),  \ \ \ \ \ \ \hat B\circ \hat D
\sim  z^rR^n(z) \circ z^n,$$
where $R(z)$ is a polynomial, $r\geq 0,$ $n\geq 1,$ and
$\GCD(n,r)=1,$ or $$\hat A\circ \hat C \sim T_n \circ T_m, \ \ \ \ \ \ \hat B\circ \hat D\sim T_m \circ T_n,$$
where $T_n,T_m$ are the corresponding Chebyshev polynomials, $n,m\geq 1,$ and \linebreak $\GCD(n,m)=1.$
Furthermore, the second
Ritt theorem remains true for arbitrary polynomial solutions of \eqref{-0}
if to replace equalities $\deg \mu_1=1, \deg \mu_2=1$ by
the equalities $$\deg \mu_1=\GCD(\deg A,\deg B),\ \ \ \deg \mu_2=\GCD(\deg C,\deg D)$$
(see \cite{en}, \cite{tor}).

Notice that the classification of polynomial solutions of \eqref{-0} appears in a variety of different
contexts some of which are quite unexpected. For example, this classification is closely related to the problem of description of
Diophantine equations of the form $A(x)=B(y),$ $A,B\in \mathbb Z[z],$ having an infinite number of integer solutions (see \cite{f1}, \cite{bilu}), and to the
problem of description of
polynomials $C,D$ satisfying the equality $C^{-1}\{S\}=D^{-1}\{T\}$ for some compact sets $S, T\subset \C$, recently solved in \cite{p1}. Notice also that the problem of description of
solutions of \eqref{-0} such that $C$ and $D$ are polynomials while $A,B$ are allowed to be arbitrary
rational (or even just continuous)
functions on the sphere can be reduced to the description
of polynomial solutions
(see \cite{p2}).
A more detailed account of different results related to the second Ritt theorem can be found
in the recent papers \cite{pppp}, \cite{pak}.

The classification of polynomial solutions of \eqref{-0} essentially reduces to the description of polynomials $A,B$
such that the algebraic curve \be \la{cu} A(x)-B(y)=0\ee
has an irreducible factor of genus zero with one point at infinity.
On the other hand, the proof of the first Ritt theorem
can be given in purely algebraic terms which do not involve the genus
condition in any form.
Indeed, if $G(F)\leq {\rm Sym}(\Omega)$ is the monodromy group of a rational function
$F$
then equivalence classes of maximal decompositions of $F$ are in a one-to-one correspondence
with maximal chains of subgroups
\be \la{zep}
G_{\omega}(F)=T_0 < T_{1} < \ ... \ < T_r=G(F),\ee
where $G_{\omega}(F)$ is the stabilizer of an element $\omega\in \Omega$ in the group $G(F)$. Therefore,
any two maximal decompositions of $F$ are weakly equivalent if and only if
for any two maximal chains of subgroups as above ${\f R_1},$ ${\f R_2}$ there exists a collection of maximal chains of subgroups ${\f T}_i$, $1\leq i \leq s,$ such that
${\f T}_1={\f R_1},$ ${\f T}_s= {\f R_2},$ and ${\f T}_{i+1}$ is obtained from ${\f T}_i,$ $1\leq i \leq s-1,$ by a replacement of exactly one group. It was shown
in the paper \cite{mu} (Theorem R.3) that
the last condition is satisfied
for any permutation group $G$ containing an abelian transitive subgroup.
Since the monodromy group of a polynomial always contains a cyclic subgroup with one orbit (its generator
corresponds to the loop around infinity), this implies in particular the truth of the first Ritt theorem for polynomials.

It was also proved in the paper \cite{mu} (Claim 1) that if $A,B,C,D$ are indecomposable polynomials satisfying
\eqref{-0} such that the decompositions $A\circ C$ and $B\circ D$ are non-equivalent then
the groups $G(A)$ and $G(D)$ as well as the groups $G(C)$ and $G(B)$ are permutation equivalent.
Since any two maximal decompositions of a polynomial $P$
are weakly equivalent, this implies by induction that for any two maximal decompositions \eqref{decc} of $P$
there exists a permutation $\sigma\in S_k$ such that the monodromy groups of $U_i$ and $V_{\sigma(i)},$
$1\leq i \leq k,$ are permutation equivalent (\cite{mz}).
The algebraic counterpart of this fact is the following
statement: if $G\leq {\rm Sym}(\Omega)$ is a permutation group containing a cyclic subgroup with one orbit then
for any two maximal chains
$$G_{\omega}=A_0 < ... < A_k=G\ \ {\rm and}\ \ G_{\omega}=B_0 < ... < B_m=G$$
the equality $k=m$ holds and
there exists a permutation $\sigma\in S_k$ such that the
permutation group
induced by the action of $A_i$ on cosets of $A_{i-1}$ is permutation equivalent
to the
permutation group induced by the action of
$B_{\sigma(i)}$ on cosets of $B_{\sigma(i)-1},$ $1\leq i \leq k.$
If a permutation group $G$ satisfies this condition,
we will say that $G$ satisfies the
{\it Jordan-H{\" o}lder theorem for imprimitivity systems}.

In this paper, we extend the above results about the permutation groups $G$ containing a cyclic group with one orbit
to the permutation groups
containing a cyclic subgroup $H$ with at most {\it two} orbits and apply these results to
rational functions (or more generally to meromorphic functions on compact Riemann surfaces)
the monodromy group of which contains $H.$

First, we prove that
for a permutation group $G$ containing $H$ the lattice $L(G_{\omega},G)$, consisting of subgroups of $G$ containing $G_{\omega}$, is
lower semi-modular and even a stronger condition of the modularity of $L(G_{\omega},G)$ holds
whenever $L(G_{\omega},G)$ does not contain a sublattice isomorphic to the subgroup lattice of
a dihedral group.
It follows easily from the lower semimodularity of $L(G_{\omega},G)$ that one can pass from any chain of subgroups \eqref{zep} to any other such a chain by a sequence of replacements as above and therefore the first Ritt theorem extends to rational functions the monodromy group of which contains $H$. Notice that this implies in
particular that the first Ritt theorem holds
for rational functions with at most two poles. Although for such functions the result was know
previously (see \cite{ar}, \cite{pak}, \cite{zi}) the algebraic proof turns out to be more simple
and illuminating. Notice also that our description of the lattice $L(G_{\omega},G)$ for groups $G$ containing $H$
has an interesting connection with the problem of description
of algebraic curves having a factor of genus zero with at most two points at infinity,
studied in \cite{f1}, \cite{bilu}.

Further, we prove that if a permutation group $G$
contains a cyclic subgroup with two orbits of {\it different length} then the lattice $L(G_{\omega},G)$ is always modular and $G$ satisfies the Jordan-H{\" o}lder theorem for imprimitivity systems.
This implies in particular that if $F$ is a rational function which has only two poles
and the orders of these poles are different between themselves then any two maximal decompositions \eqref{decc} of $F$ have the same length and
there exists a permutation $\sigma\in S_r$ such that the monodromy groups of $U_i$ and $V_{\sigma(i)},$
$1\leq i \leq r,$ are permutation equivalent. We also show that the Jordan-H{\" o}lder theorem for imprimitivity systems holds for any permutation group containing a transitive
Hamiltonian subgroup that generalizes the corresponding results of \cite{mu}, \cite{mz}.

For arbitrary rational functions the first Ritt theorem fails to be true.
The simplest counterexamples are provided by the functions
which are regular coverings of the sphere (that is for which $G_{\omega}=e$)
with the monodromy group $A_4,$ $S_4$, or $A_5.$ These functions were described for the first time by F. Klein in \cite{klein} and nowadays can be interpreted
as Belyi functions of Platonic solids (see \cite{cg}, \cite{zv}).
For such a function
its maximal decompositions
simply correspond to maximal chains of subgroups in its monodromy group.
Therefore, since any of the groups $A_4$, $S_4,$ $A_5$ has maximal chains of subgroups of different length, for the corresponding
Klein functions the first Ritt theorem is not true.

Although the fact that
the Klein functions
provide counterexamples to the first Ritt theorem
is a well known part of the mathematical ``folklore'', the systematic description of compositional properties of these functions seems to be absent. In particular, to our best knowledge
maximal decompositions
which do not satisfy the first Ritt theorem were found explicitly
only for the Klein function corresponding to the group $A_4$ (see \cite{gs}, \cite{be2}).
In the Appendix to this paper we provide a detailed analysis of maximal decompositions
of the Klein functions and give related explicit examples of non weakly equivalent maximal decompositions. In particular,
we give an example of a rational function with {\it three} poles
having maximal decompositions of different length. This example shows that with no additional
assumptions the first Ritt theorem can not be extended to rational functions the monodromy of which
contains a cyclic subgroup with more than two orbits.

\vskip 0.2cm
\noindent{\bf Acknowledgments}. The authors would
like to thank I. Kovacs, P. M\"uller, U. Zannier, and A. Zvonkin for discussions of different questions
related to the subject of this paper.

\section{Jordan-H\"older theorem for imprimitivity systems}
\subsection{Lattices, imprimitivity systems, and decompositions of functions}

Recall that {\it a lattice} is a partially ordered set \cmnt{ $(L,\leq)$ }in which every pair of elements $x,y$
has a unique supremum $x\vee y$ and an infimum $x\wedge y$ (see e.g. \cite{ai}).
Our basic example of a lattice is a lattice $L(G)$ of all subgroups of a group $G,$ where
by definition $G_1\leq G_2$ if $G_1$ is a subgroup of $G_2$ (clearly, $G_1\cap\, G_2$ is an infimum of $G_1,G_2$ and
$\langle G_1,G_2\rangle$ is a supremum). A simplest example of the lattice $L(G)$ is obtained if
$G$ is a cyclic group of order $n$. In this case $L(G)$ is isomorphic
to the lattice $L_n$ consisting of
all divisors of $n$, where by definition $d_1\leq d_2$ if $d_1\vert d_2.$

A {\it sublattice} of a lattice $L$ is a non-empty subset $M\subseteq L$
closed with respect to $\vee$ and $\wedge$. For example, for any subgroup $H$ of a group $G$
the set $$L(H,G):=\{X\leq G\,|\,H\leq X\leq G\}$$ is a sublattice of $L(G)$.
Another example of a sublattice of $L(G)$ is the lattice $$L(A,AB):=\{X\leq G\,|\, A\leq X\subseteq AB\}$$
(notice that in our notation $X\leq G$ means that $X$ is a subgroup of
$G$ while $X\subseteq AB$ means that $X$ is a subset of the set $AB$ which in general is not supposed to be a group). Recall that by the Dedekind identity (see e.g. \cite{huppert}, p. 8)
for arbitrary subgroups $A,B,X$ of a group $G$ such that $A\leq X\subseteq AB$
the equality $X=A(X\cap B)$ holds. It follows from the Dedekind identity that the mapping $f:X\mapsto X\cap B$
is a monomorphism from the lattice $L(A,AB)$ into
the lattice $L(A\cap B,B)$ with the image consisting of all subgroups of $B$ which are permutable
with $A$. We will call $f$ the Dedekind monomorphism.

For elements $a,b$ of a lattice $L$ the symbol $a<\cdot\ b$ denotes that
$a\leq b$
and there exists no element $c\neq a,b$ of $L$ such
that $a\leq c \leq b.$
A lattice $L$ is called {\it semimodular} \cmnt{\cite{ai}} if for any $a,b\in L$ the condition
\be \la{m1} a\wedge b<\cdot\ a, \ \ \ \ \ a\wedge b<\cdot\ b,\ee
imply the condition \be \la{m2} b <\cdot\ a\vee b, \ \ \ \ \ a <\cdot\ a\vee b.\ee
If vice versa condition \eqref{m2} implies condition \eqref{m1}, the lattice $L$ is called {\it lower semimodular}.
A lattice $L$ is called {\it modular} if $L$ is semimodular and lower semimodular.
{\it A maximal chain} $\f R$ between elements $a,b$ of $L$ is a collection
$a_0,a_2,\dots a_k$ of elements of $L$ such that
$$\f R: \  a=a_0<\cdot\ a_1<\cdot\ \ \ \dots \ \ <\cdot\ a_k=b.$$
The number $k$ is called the length of the chain $\f R$
(we always assume that in the lattices considered the length of a chain between $a$ and $b$
is uniformly bounded by a number depending on $a$ and $b$ only).

It is well known (see e.g. \cite{ai}) that for a semimodular or lower semimodular lattice all maximal chains between
two elements have the same length. Below, using essentially the same proof,
we give a modification of this statement in the spirit of the first Ritt theorem.

Say that two maximal chains between elements $a$ and $b$ of a lattice $L$ are $r$-equivalent
if there exists a sequence of maximal chains
$\f T_1, \f T_2,\, \dots \, \f T_s$ between $a,b$ such that $\f T_1=\f R_1,$ $\f T_s=\f R_2,$ and
$\f T_{i+1}$ is obtained from $\f T_{i},$ $1\leq i \leq s-1,$ by a
replacement of exactly one element. Clearly, all $r$-equivalent chains
have an equal length.

\bt \la{ritt} Let $L$ be a semimodular or lower semimodular lattice.
Then any two maximal chains between any elements $a$ and $b$ of $L$ are $r$-equivalent.
\et
\pr Since after the inversion of the ordering of a lattice the condition of semi-modularity
transforms to the condition of lower semi-modularity and vice versa, it
is enough to prove the theorem for lower semi-modular lattices.

Fix $a\in L$. For arbitrary $b\in L$ denote by $d(b)$ a maximum of lengths of maximal chains between
$a$ and $b$. We will prove the theorem by induction on $d(b).$
For $b$ satisfying $d(b)\leq 1$ the theorem is obviously true. Suppose that the theorem is proved for $b$ satisfying $d(b)\leq n-1$ and let
$$\f R_1: \  a=a_0<\cdot\ a_2<\cdot\ \ \ \dots \ \ <\cdot\ a_{k_1}=b, \ \ \ \f R_2: \  a=b_0<\cdot\ b_2<\cdot\ \ \ \dots \ \ <\cdot\ b_{k_2}=b$$
be two maximal chains between $a$ and an element $b\in L$ such that $d(b)=n.$

If $a_{k_1-1} = b_{k_2-1}$, then we are done by induction. So, we may assume that $a_{k_1-1} \neq b_{k_2-1}$.
Then by the maximality of $a_{k_1-1}$ and $b_{k_2-1}$ in $b$ we conclude $a_{k_1-1} \vee b_{k_2-1}=b$.
Hence
$$a_{k_1-1} <\cdot\ a_{k_1-1} \vee b_{k_2-1}, \ \ \ b_{k_2-1} <\cdot\ a_{k_1-1} \vee b_{k_2-1}$$
and therefore by the lower semi-modularity of $L$ we have:
\be \la{as} a_{k_1-1} \wedge b_{k_2-1} <\cdot\ a_{k_1-1}, \ \ \  a_{k_1-1} \wedge b_{k_2-1}<\cdot\ b_{k_2-1}.\ee
Let $$a=c_0<\cdot\ c_2<\cdot\ \ \ \dots \ \ <\cdot\ c_{l}=a_{k_1-1} \wedge b_{k_2-1}$$
be any maximal chain between $a$ and $a_{k_1-1} \wedge b_{k_2-1}$ and
\be\la{ikk}  a=c_0<\cdot\ c_2<\cdot\ \ \ \dots \ \ <\cdot\ c_{l}<\cdot\ a_{k_1-1}\ee be
its extension to a maximal chain between $a$ and $a_{k_1-1}$. Since $d(a_{k_1-1})$ is obviously
less than $d(b)$, it follows from the induction
assumption that the chain
$$a=a_0<\cdot\ a_2<\cdot\ \ \ \dots \ \ <\cdot\ a_{k_1-1}$$
obtained from $\f R_1$ by deleting $a_{k_1}$ is $r$-equivalent
to the chain \eqref{ikk}.
Therefore, the chain $\f R_1$ and the chain
\be \la{y1} a=c_0<\cdot\ c_2<\cdot\ \ \ \dots \ \ <\cdot\ c_{l}<\cdot\ a_{k_1-1}<\cdot\ b \ee
also are $r$-equivalent.

Similarly, the chain
$\f R_2$ is $r$-equivalent to the chain
\be \la{y2} a=c_0<\cdot\ c_2<\cdot\ \ \ \dots \ \ <\cdot\ c_{l}<\cdot\ b_{k_2-1}<\cdot\ b.\ee
Since chains \eqref{y1} and \eqref{y2} are $r$-equivalent, we conclude that the chain $\f R_1$ is $r$-equivalent to the chain $\f R_2$. \qed

\br Notice that there exist lattices which are not semimodular or lower semimodular
such that any two maximal chains between any elements are $r$-equivalent. An example of such a
lattice is shown on Fig. \ref{la}. \vskip 0.2cm

\begin{figure}[t]
\renewcommand{\captionlabeldelim}{.}
\medskip
\epsfxsize=1.8truecm
\centerline{\epsffile{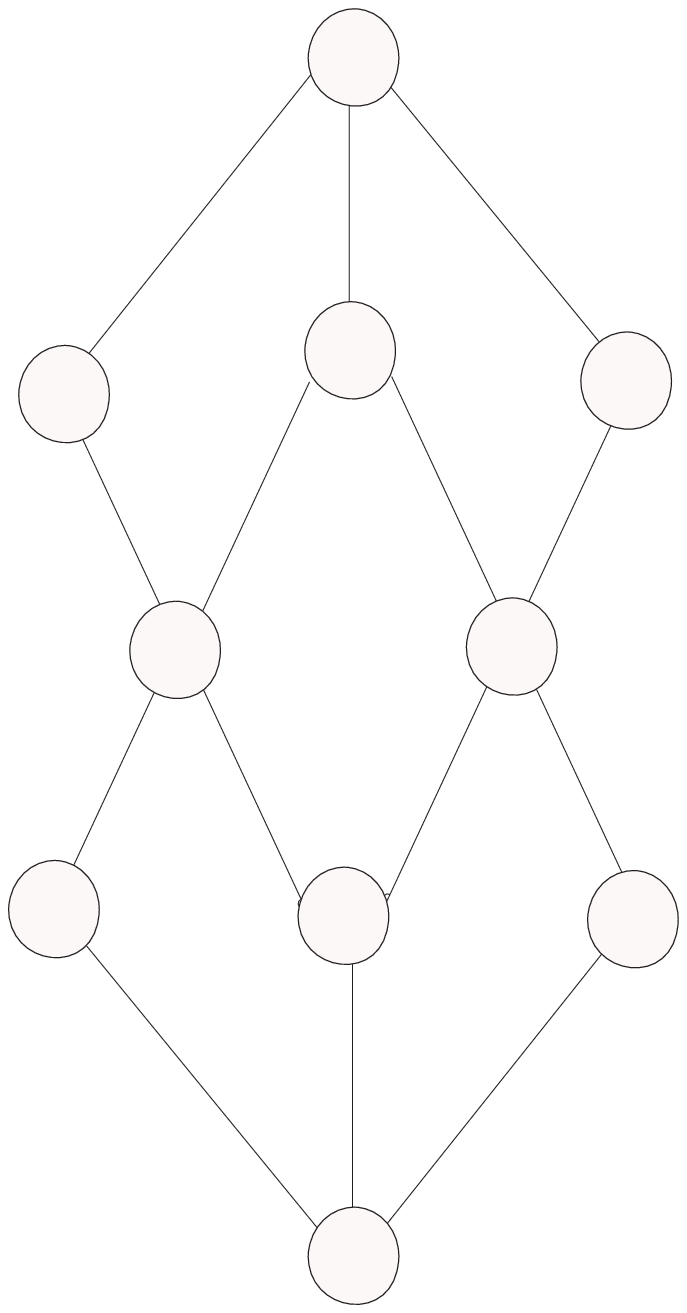}}
\caption{}\la{la}
\medskip
\end{figure}

Let $\Omega$ be a finite set and $G\leq \sym(\Omega)$ be a transitive permutation
group. Recall that a partition
$\f E$ of $\Omega$ is called an \New{imprimitivity system} of $G$ if $\f E$ is $G$-invariant.
Elements of $\f E$ are called \New{blocks}. For a point $\omega\in\Omega$ we will denote by ${\f E}(\omega)$ a unique block of $\f E$ which
contains $\omega.$ Since the group $G$ permutes the elements
of $\f E$ transitively, all blocks of $\f E$ have the same cardinality denoted by $n_{{\f E}}$.
Denote
by ${\f E}(G)$ the set of all imprimitivity systems of $G$. It is a
partially ordered set, where by definition ${\f E}\leq {\f F}$
if ${\f E}$ is a refinement of ${\f F}$. Notice that if ${\f E}\leq {\f F}$
then $n_{{\f F}}/n_{{\f E}}$ is an integer denoted by $[{\f F}:{\f E}]$.

It is easy to see that
${\f E}(G)$ is a lattice where the lattice operations are defined as follows
$$
{\f E}\land {\f F}:=\{\Delta\cap\Gamma\,|\,\Delta\in{\f E},\Gamma\in{\f F}\mbox{ and }
\Delta\cap\Gamma\neq\emptyset\},
$$
$$
{\f E}\lor {\f F}:= \bigwedge\{{\f D}\in{\f E}(G)\,|\, {\f E}\leq {\f D}\mbox{ and }
{\f F}\leq{\f D}\}.
$$
It is well known that the lattice ${\f E}(G)$ is isomorphic to the subgroup lattice $L(G_\omega, G)$ where
$\omega\in\Omega$ is an arbitrary fixed point. The correspondence between two
sets is given by the formula ${\f E}\mapsto G_{{\f E}(\omega)}$, where
$$G_{{\f E}(\omega)}:=\{g\in G\,|\, {\f E}(\omega)^g={\f E}(\omega)\}.$$
Vice versa, an imprimitivity system corresponding to a
subgroup $K\in L(G_\omega,G)$ is defined as follows ${\f E}_K:=\{\omega^{Kg}\,|\,g\in G\}$.
Notice that for any ${\f E},{\f F}\in {\f E}(G)$ we have:
$$
G_{({\f E}\land{\f F})(\omega)} = G_{{\f E}(\omega)} \cap G_{{\f F}(\omega)}, \ \
G_{({\f E}\lor {\f F})(\omega)} = \sg{G_{{\f E}(\omega)},G_{{\f F}(\omega)}}
.$$ Moreover, if
${\f E}\leq {\f F}$ then
$[{\f F}:{\f E}]=[G_{{\f F}(\omega)}:G_{{\f E}(\omega)}].$

If a group $G$ is the monodromy group of a rational function $F$, then imprimitivity systems of $G$
are in a one-to-one correspondence with
equivalence classes of decompositions $A\circ B$ of $F$. Namely, suppose that
$G$ is realized as a permutation group acting
on the set $z_1,z_2,\dots,z_n$ of preimages of a non critical value $z_0$ of $F=A\circ B$
under the map $F\,:\C\P^1\rightarrow \C\P^1$, and let $x_1,x_2,\dots,x_r$ be the set of preimages of $z_0$ under the map $A\,:\C\P^1\rightarrow \C\P^1$.
Then blocks of the imprimitivity system of $G$ corresponding to the equivalence class of  decompositions of $F$ containing $A\circ B$, are just preimages of
the points $x_1,x_2,\dots,x_r$ under the map $B\,:\C\P^1\rightarrow \C\P^1$.
More generally, equivalence classes of decompositions of a rational function $F$ are in a one-to-one correspondence with chains
of subgroups
$$
G_{\omega}=T_0< T_{1} < \ ... \ < T_r=G,$$
where $G$ is the
monodromy group of $F$.

Following \cite{pak} we say that two maximal decompositions $\f D_1,\f D_2$ of a rational function $F$ are weakly equivalent
if there exists
a chain of maximal decompositions $\f F_i$, $1\leq i \leq s,$ of $F$ such that
$\f F_1=\f D_1,$ $\f F_s\sim \f D_2$, and $\f F_{i+1}$ is obtained from $\f F_i,$
$1\leq i \leq s-1,$
by replacing two successive functions in $\f F_i$ by other functions with the same composition.
The remarks above imply that two maximal decompositions of
$F$ are weakly equivalent if and only if corresponding maximal chains in $L(G_{\omega},G)$
are $r$-equivalent.
In particular, the conclusion of the first Ritt theorem is true for a rational function $F$
if and only if all maximal chains between $G_{\omega}$ and $G$ in $L(G_{\omega},G)$ are
$r$-equivalent.
Therefore, Theorem \ref{ritt} implies the following corollary (cf. \cite{pak}, Th. 2.5).

\bc \la{df} Let $F$ be a rational function  such that the lattice $L(G_\omega, G),$ where $G$ is the
monodromy group of $F$, is
semi-modular or lower semi-modular. Then all maximal decompositions of $F$ are weakly equivalent. \qed
\ec

The Corollary \ref{df} shows that the groups $G$ for which $L(G_{\omega},G)$ is semi-modular or lower semi-modular
are of special interest for factorization theory
of rational functions. The simplest examples of such groups are groups containing a transitive cyclic subgroup.

\bt \la{cyc} Let $G\leq S_n$ be a permutation group containing a transitive cyclic subgroup
$C_n$. Then the lattice
$L(G_{1},G)$ is a modular lattice isomorphic to a sublattice of the lattice $L_n.$
\et

\pr Since any sublattice of a modular lattice is modular (see e.g. \cite{ai}) and
it is easy to see that $L_n$ is
modular, it is enough to prove that $L(G_{1},G)$ is isomorphic to a sublattice of $L_n.$

The transitivity of $C_n$ implies that $G=G_{1} C_n$.
Therefore, the Dedekind monomorphism $f:X\mapsto X\cap C_n$ maps
$L(G_{1},G)$ into a sublattice of $L(G_{1}\cap C_n,C_n)$. On the other hand,
$$L(G_{1}\cap C_n,C_n)=L(e,C_n) \cong L_n.\ \ \Box$$
\vskip 0.2cm
Note that Theorem \ref{cyc} implies the following proposition
(cf. \cite{en}, \cite{tor}).

\bc \la{eng}
Let $A,B,C,D$ be polynomials such that $$A\circ C=B\circ D.$$ Then there exist polynomials
$U, V, \hat A, \hat C, \hat B, \hat D, $  where
$$\deg U=\GCD(\deg A,\deg B),  \ \ \ \deg V=\GCD(\deg C,\deg D),$$
such that
$$A=U\circ \hat A, \ \  B=U\circ \hat B, \ \ C=\hat C\circ V, \ \  D=\hat D\circ V,$$
and $$\hat A\circ \hat C=\hat B\circ \hat D.$$ In particular, if $\deg A=\deg B$ then the decompositions
$A\circ C$ and $B\circ D$ are necessarily equivalent.
\ec

\subsection{Jordan-H\"older theorem for groups with normal imprimitivity systems}
Let as above $G$ be a transitive permutation group. It is easy to see that if
$N$ is a normal subgroup of $G$ then its orbits form an imprimitivity
system of $G$. Such an imprimitivity system is called \New{normal} and is denoted by $\Omega/N$.
For an imprimitivity system ${\f E}\in{\f E}(G)$
set
$$G_{{\f E}}:=\{g\in G\,|\,\forall_{\Delta\in{\f E}}\ \Delta^g=\Delta\}.$$
Notice that each block of ${\f E}$ is a union of $G_{{\f E}}$-orbits and
$G_{{\f E}} = \core_G(G_{{\f E}(\omega)}).$
In particular, $G_{{\f E}}$ is a normal subgroup of $G$.

Let us call a subgroup $A\in L(G_\omega,G)$ \New{core-complementary}
if $A = G_\omega \core_G(A)$.

\bp\la{normality}
An imprimitivity system ${\f E}\in{\f E}(G)$ is normal if and only if the group $G_{{\f E}(\omega)}$
is core-complementary.
\ep
\pr
Indeed, if
\be \la{kott} G_{{{\f E}(\omega)}}=G_\omega  \core_G(G_{{{\f E}(\omega)}})=
G_\omega G_{{\f E}} \ee
then
$${{\f E}(\omega)} = \omega^{G_{{\f E}(\omega)}} = \omega^{G_{{\f E}}}$$
and hence $G_{{\f E}}$ acts transitively on ${{\f E}(\omega)}$.
Since $G_{{\f E}}\normaleq G$, this implies that $G_{{\f E}}$ acts transitively on every block of ${\f E}$. Thus blocks
of ${\f E}$ are orbits of the normal subgroup $G_{{\f E}}$.

Vice versa, if ${\f E}$ is normal then ${\f E} = \Omega/N$ for some $N\normaleq G$.
This implies that $G_{{\f E}(\omega)} = G_\omega N$ and $N\leq G_{{\f E}}$. It follows now from
$$
G_{{\f E}(\omega)} = G_\omega N \leq G_\omega G_{{\f E}} \leq G_{{\f E}(\omega)}
$$
that equality \eqref{kott} holds.
\qed

Recall that two subgroups $A$ and $B$ are called \New{permutable} if $AB=BA$, or, equivalently,
$\sg{A,B} = AB$. Recall also that if $A$ and $B$ are subgroups of finite index of $G$  then the inequality
\begin{equation}\la{form}
\left[\sg{A,B}:B\right]\geq\left[A:A\cap B\right]
\end{equation}
holds and the equality in \eqref{form} attains if and only if $A,B$ are permutable (see e.g. \cite{kur}, p. 79).

Denote by $L_c(G_\omega,G)$ the subset
of $L(G_\omega,G)$ consisting
of all core-complemen\-tary subgroups. Notice that in general $L_c(G_\omega,G)$ is {\it not} a
sublattice of $L(G_\omega,G)$

\bp\la{properties} The following conditions hold:

\begin{enumerate}
\item[(a)] If $A\in L_c(G_\omega,G)$, then $AB=BA$ for each $B\in L(G_\omega,G)$;

\item[(b)] If $A,B\in L_c(G_\omega,G)$, then $AB\in L_c(G_\omega,G)$.
\end{enumerate}

\ep

\pr (a) In order to lighten the notation set $N=\core_G(A)$. In view of Proposition \ref{normality} we have: $$AB = G_\omega N B =
N G_\omega B = NB=B N = B G_\omega N = BA.
$$
(b) Set $M=\core_G(B)$. Since $MN\normaleq G$ and $MN\leq AB$, we have: $$MN\leq\core_G(AB).$$
It follows now from Proposition \ref{normality} that
$$AB = G_\omega N G_\omega M= G_\omega MN\leq G_\omega \core_G(AB)\leq AB.$$ Therefore,
$G_\omega \core_G(AB)= AB$ and hence
$AB\in L_c(G_\omega, G)$ by Proposition \ref{normality}.
\qed

\bp \la{lp} Let $A,B\leq G$ be permutable subgroups. If $A\cap B$ is maximal in
$A$ and $B$, then $A$ and $B$ are maximal in $\sg{A,B} = AB$.
\ep

\pr Let $A_1$ be a subgroup of $G$ satisfying $A \leq A_1 \leq AB$. It follows from
$$A\cap B\leq A_1\cap B \leq B$$ that either $A_1\cap B = A\cap B$
or $A_1\cap B= B$. It follows now from the
Dedekind identity $A_1 = A (A_1\cap B)$ that in the first case $A_1=A$ while in the second one $A_1 = AB$.\qed

\bp \la{modular} If any two subgroups of $L(G_\omega,G)$ are permutable,
then the lattice $L(G_\omega,G)$ is modular.
\ep
\pr Indeed, if $A\cap B$ is maximal in $A$ and $B$
then $A$ and $B$ are maximal in $\sg{A,B} = AB$
by Proposition~\ref{lp}.

Suppose now that $A$ and $B$ are maximal
in $AB$ and let $A_1$ be a subgroup of $G$ satisfying $A\cap B\leq A_1\leq A$. Then
$$B\leq A_1B\leq AB$$ implies that
either $B=A_1B$ or $A_1B = AB$. If $B=A_1B$,
then $A_1\leq B$ and therefore $A_1 = A\cap B.$
On the other hand, if $A_1B = AB$ then it follows from $A\leq AB=A_1B$ that for any
$a\in A$ there exist $a_1\in A_1$ and $b\in B$ such that $a=a_1 b$. Since the last equality yields that
$b\in A\cap B$, this implies that $A\leq A_1(A\cap B)\leq A_1$ and hence $A_1=A.$ \qed

Let $H\leq G$ be an arbitrary subgroup and $H\backslash G:=\{Hx\,|\,x\in G\}$.
Denote by $G//H$ a permutation group arising from the natural action of $G$ on
$H\backslash G$. Thus $G//H$ is always considered as a subgroup of $\sym(H\backslash G)$.
Notice that if $N\normaleq G$ is contained in $H$, then the groups $G//H$ and $(G/N)//(H/N)$
are permutation equivalent.
Below we will denote permutation equivalence by $\cong_p$.

Say that a transitive permutation group $G\leq\sym(\Omega)$ satisfies the
Jordan-H{\" older} theorem for imprimitivity systems if any two maximal chains
$$G_\omega=A_0 < ... < A_k=G\ \ {\rm and} \  \ G_\omega=B_0 < ... < B_m=G$$
of the lattice $L(G_\omega,G)$ have the same length (that is $k=m$)
and
there exists a permutation $\sigma\in S_k$ such that the
permutation groups
$A_i//A_{i-1}$ and $B_{\sigma(i)}//B_{\sigma(i)-1},$ $1\leq i \leq k,$ are permutation equivalent. Notice that if $G$ is the monodromy group of a rational function $F$ then
it follows from the
correspondence between imprimitivity systems of $G$ and
equivalence classes of decompositions of $F$ that $G$ satisfies the
Jordan-H{\" older} theorem for imprimitivity systems if and only if
any two maximal
decompositions of $F$
$$F=U_1\circ U_2\circ \dots \circ U_k \ \ {\rm and} \  \ F=V_1\circ V_{2}\circ \dots \circ V_m,$$
have the same length and the there exists a permutation $\sigma\in S_k$ such that the monodromy groups of $U_i$ and $V_{\sigma(i)},$
$1\leq i \leq k,$ are permutation equivalent.

\bt\la{JH} Let $G$ be a permutation group such that $L(G_\omega,G) = L_c(G_\omega,G).$ Then the lattice $L(G_\omega,G)$ is modular and
$G$ satisfies the Jordan-H{\" o}lder
theorem for imprimitivity systems.
\et
\pr First of all observe that since by Proposition~\ref{properties} any two subgroups of $L(G_\omega,G)$ are
permutable it follows from Proposition~\ref{modular} that $L(G_\omega,G)$ is a modular lattice.
Let now
$${\f A}:=G_\omega=A_0 < ... < A_k=G,\ \ {\rm and} \ \ {\f B}:=G_\omega=B_0 < ... < B_m=G$$
be two maximal chains of $L(G_\omega,G)$.
Since $L(G_\omega, G)$ is a modular lattice, it follows from Theorem~\ref{ritt}
that $k=m$ and $\f A$ and $\f B$ are $r$-equivalent.
Therefore by induction it is sufficient to prove the theorem for the case when ${\f B}$ and ${\f A}$
differs at exactly one place, say $i$ ($1\leq i < k$). Clearly, in this case we have:
$$ A_{i-1}=B_{i-1}=A_i\cap B_i, \ \ \ A_{i+1} = B_{i+1} = A_iB_i.$$

In order to lighten the notation set $$N:=\core_{A_{i+1}}(A_i).$$ It follows from the equality $A_i = G_\omega\core_G(A_i)$ that $A_i = A_{i-1}\core_G(A_i)$. Therefore,
$$A_i = A_{i-1}\core_G(A_i)\leq A_{i-1}N\leq A_i$$ and hence
$$A_i =  A_{i-1}N=B_{i-1} N$$ and
$$A_{i+1} = A_iB_i=A_{i-1} NB_i= B_{i-1}N B_i=B_iN.$$
Since $N\leq A_i=B_{i-1} N$ and $N\normaleq A_{i+1} = B_iN,$ this implies that
$$
A_{i+1}//A_i = (B_i N)//(B_{i-1}N)\cong_p (B_i N)/N\,//\,(B_{i-1} N)/N.
$$

By the Second Isomorphism Theorem
the group $(B_i N)/N$ is isomorphic to the group $B_i/(B_i\cap N)$ and the image of $(B_{i-1} N)/N$
under this isomorphism is $$B_{i-1}(B_i\cap N)/(B_i\cap N).$$
Furthermore, it follows from $N\leq A_i$
that $B_i\cap N \leq A_i\cap B_i = B_{i-1}$. Therefore,
$$B_{i-1}(B_i\cap N)/(B_i\cap N) = B_{i-1}/(B_i\cap N)$$ and hence
$$(B_i N)/N\,//\,(B_{i-1} N)/N\cong_p B_i/(B_i\cap N)\,//\,B_{i-1}/(B_i\cap N).$$
Finally, since $B_i\cap N\normaleq B_{i}$,
$$B_i/(B_i\cap N)\,//\,B_{i-1}/(B_i\cap N)\cong_p B_i//B_{i-1}$$ and hence
$
A_{i+1}//A_i \cong_p B_i//B_{i-1}.$
Replacing $A$ and $B$ in the above argument
we obtain similarly that $B_{i+1}//B_i \cong_p A_i//A_{i-1}$.
\qed

Recall that a group is called \New{Hamiltonian} if all its subgroups are normal.

\bt\la{hamiltonian} Let $G$ be a permutation group containing a transitive Hamiltonian subgroup $K$.
Then $L(G_\omega,G)$ is a modular lattice isomorphic to a sublattice of $L(K)$ and
$G$ satisfies the Jordan-H{\" o}lder theorem for imprimitivity systems.
\et
\pr
It follows from the transitivity of $K$ that $G=G_\omega K$. Therefore, by the Dedekind
monomorphism
$L(G_\omega,G)$ is isomorphic to a sublattice of $L(G_\omega\cap K,K).$ Clearly,
$G_\omega\cap K=K_\omega.$ Furthermore, since $K$ is Hamiltonian,
the subgroup $K_\omega$ is normal in $K$ and therefore for
any $\omega^{\prime}\in \Omega$ the equality $K_\omega=K_{\omega^{\prime}}$ holds. This implies that $K_\omega = 1$ and hence
$L(G_\omega, G)$ is isomorphic to a sublattice of $L(K)$. Since $L(K)$ is modular by Proposition \ref{modular},
$L(G_\omega, G)$ is modular as well.

By Theorem~\ref{JH} in order to prove that $G$ satisfies the Jordan-H{\" o}lder theorem it is enough to show that
$L_c(G_\omega,G) = L(G_\omega, G)$.
Observe first that it follows from
$G = G_\omega K$ that for arbitrary $A\in L(G_\omega,G)$ the equality \be \la{wwee} \core_G(A)=\cap_{g\in K} gAg^{-1}\ee holds.
On the other hand, since $K$ is Hamiltonian, $A\cap K\normaleq K$. Therefore,
for each $g\in K$ we have:
$$g^{-1}(A\cap K)g =A\cap K\leq A$$ implying
\be \la{eeww} A\cap K \leq gAg^{-1}.\ee It follows now from \eqref{wwee} and \eqref{eeww} that $A\cap K\leq \core_G(A)$ and hence
$$G_\omega(A\cap K)\leq G_\omega\core_G(A)\leq A.$$ Since by Dedekind's identity $G_\omega(A\cap K)=A$, we conclude that $G_\omega\core_G(A)=A$ for any $A\in L(G_\omega,G)$ and hence $L_c(G_\omega,G) = L(G_\omega, G)$ by Proposition \ref{normality}.
\qed

\bc \la{ccoo} Let $F$ be a rational function such that its monodromy group contains a transitive Hamiltonian subgroup. Then any two maximal decompositions of $F$ are weakly equivalent. Furthermore, for any two decompositions of $F$:
$$F=U_1\circ U_2\circ \dots \circ U_k \ \ {\it and} \  \ F=V_1\circ V_{2}\circ \dots \circ V_k,$$
there exists a permutation $\sigma\in S_k$ such that the monodromy groups of $U_i$ and $V_{\sigma(i)},$
$1\leq i \leq k,$ are permutation equivalent.
\ec

Notice that the condition of
Corollary \ref{ccoo} is satisfied in particular if $K$ is cyclic or abelian. Therefore, Corollary \ref{ccoo}
generalizes Theorem R.3 and Claim 1 of
\cite{mu}, and Theorem 1.3 of \cite{mz}.
\vskip 0.2cm

\subsection{Jordan-H\"older theorem for groups containing a cyclic subgroup with two orbits of different length}
Let $\Omega$ be a finite set, $h\in \sym(\Omega)$ be a permutation which is a product of exactly
two disjointed cycles, and $H:=\sg{h}$. For the rest of this subsection it is assumed that $G\leq\sym(\Omega)$ is a transitive
permutation group containing $H$.
Without loss of generality we may assume that $G\leq S_n$ and
$$h=(1\,2\,\dots\, n_1)(n_1+1\, n_1+2\, \dots\, n_1+n_2),$$ where
$1\leq  n_1,n_2 < n,$ $n_1+n_2=n.$

Say that an imprimitivity system ${\f E}\in{\f E}(G)$
is $H$-transitive (resp. $H$-intransitive) if the action of $H$ on blocks of ${\f E}$ is
transitive (resp. intransitive). Say that a group $K\in L(G_{\omega},G)$ is $H$-transitive (resp. $H$-intransitive) if the corresponding ${\f E_K}\in{\f E}(G)$ is $H$-transitive (resp. $H$-intransitive).

Since $H$ permutes blocks of $\f E$, it is easy to see that if
${\f E}$ is $H$-transitive
then there exist numbers $d\vert n$ and $i_1,i_2$ $1\leq i_1,i_2\leq d,$ such that
any block of $\f E$ is equal to $W_{i_1,d}^1\cup W_{i_2,d}^2,$ where the symbol
$W_{j,l}^1$ (resp. $W_{j,l}^2$) denotes a union of numbers from the
segment $[1,n_1]$ (resp. from the segment $[n_1+1,n_1+n_2]$) equal to $j$
by modulo $l.$ On the other hand, if
${\f E}\in{\f E}(G)$ is $H$-intransitive then there exist numbers
$d_1\vert n,d_2\vert n$ and $i_1,i_2,$ $1\leq i_1\leq d_1,$ $1\leq i_2\leq d_2,$
such that
\be \la{rav} n_1/d_1=n_2/d_2=n_{\f E} \ee
and
any block of $\f E$ is equal either to $W_{i_1,d_1}^1$ or to
$W_{i_2,d_2}^2$.

\bp\la{normal} Any $H$-intransitive imprimitivity system ${\f E}\in{\f E}(G)$ is normal.
\ep
\pr In the notation above set $r={\rm lcm}(d_1,d_2)$ and $K:=\sg{h^r}$. Clearly, we have $K\leq G_{{\f E}}$ and therefore any orbit of $G_{{\f E}}$ is a union of orbits of $K.$
The length of any orbit of $K$ on $[1,n_1]$ is equal to $$\frac{n_1}{{\rm gcd}(n_1,r)}=\frac{n_{\f E}}
{\gcd(n_{\f E},r/d_1)}.$$ On the other hand,
the length of any orbit of $K$ on $[n_1+1,n_1+n_2]$ is equal to
$$\frac{n_2}{{\rm gcd}(n_2,r)}=\frac{n_{\f E}}
{\gcd(n_{\f E},r/d_2)}.$$
Therefore,
the length of any orbit of $G_{{\f E}}$ on $\Omega$ is divisible by
$${\rm lcm}\left(\frac{n_{\f E}}
{\gcd(n_{\f E},r/d_1)},\frac{n_{\f E}}
{\gcd(n_{\f E},r/d_2)} \right)=\frac{n_{\f E}}{\gcd(n_{\f E},\gcd(r/d_1,r/d_2))}=n_{\f E}.$$
This implies that orbits of $G_{{\f E}}$ coincide with blocks of ${\f E}$
and hence $\f E$ is normal. \qed

\bp\la{normal1} If $H$-transitive imprimitivity system ${\f E}\in{\f E}(G)$ is not normal,
then $n_1=n_2$ and there exists a
normal imprimitivity system ${\f E}'\leq {\f E}$ such
that $[{\f E}:{\f E}']=2$. Furthermore, ${\f E}'$ is $H$-intransitive, its blocks coincide
with the orbits of $G_{{\f E}}$,
and for any $H$-intransitive imprimitivity system ${\f F}\in{\f E}(G)$
such that ${\f F}\leq {\f E}$ we have ${\f F}\leq {\f E}'$.
\ep
\pr In the notation above set $K=\sg{h^d}.$ Clearly, any block
$W_{i_1,d}^1\cup W_{i_2,d}^2$ of ${\f E}$ is a union of exactly two orbits of $K$ and $K\leq G_{{\f E}}.$
Since ${\f E}$ is not normal,
this implies that orbits of $G_{{\f E}}$ coincide with orbits of $K$.
In particular, since orbits of $G_{{\f E}}$ have the same length the same is true for orbits of $K$
and hence $n_1=n_2.$ The rest statements of the proposition are now obvious. \qed

\bt\la{JH2} If a transitive permutation group $G$ contains a cyclic subgroup with two orbits of different length, then $L(G_\omega,G)$ is modular and
$G$ satisfies the Jordan-H{\" o}lder theorem for imprimitivity systems.
\et

\pr It follows from Propositions \ref{normal} and \ref{normal1} that $L(G_\omega,G) = L_c(G_\omega,G)$.
Now the theorem follows from Theorem~\ref{JH}.
\qed

\bc \la{ccoo1} Let $F$ be a rational function such that $F$
has only two poles and the orders of these poles are different between themselves.
Then any two maximal decompositions of $F$ are weakly equivalent. Furthermore, for any two decompositions of $F$:
$$F=U_1\circ U_2\circ \dots \circ U_k \ \ {\it and} \  \ F=V_1\circ V_{2}\circ \dots \circ V_k,$$
there exists a permutation $\sigma\in S_k$ such that the monodromy groups of $U_i$ and $V_{\sigma(i)},$
$1\leq i \leq k,$ are permutation equivalent.
\ec

\section{The lattice of imprimitivity systems for groups containing a cyclic subgroup with two orbits}
\subsection{Semimodularity and modularity of $L(G_\omega,G)$}

\bp \la{jk} Let $G$ be a transitive permutation group. Suppose that $L(G_\omega,G)$ contains subgroups $E, F$
such that $[E:E\cap F] = [F:E\cap F]=2.$ Then $E\cap F$ is normal in $\sg{E,F}$ and $\sg{E,F}/E\cap F\cong D_{2m},$ where $2m:=[\sg{E,F}:E\cap F]$.
Furthermore, $L(E\cap F,\sg{E,F})\cong L(D_{2m}).$

\ep
\pr
Since $[E:E\cap F] = [F:E\cap F]=2$, the subgroup
$E\cap F$ is normal in $E$ and $F$ simultaneously and
therefore $E\cap F\normaleq \langle E,F\rangle.$
Since $$\langle E,F\rangle/(E\cap F) = \sg{E/(E\cap F), F/(E\cap F)}$$ and
$E/(E\cap F)\cong \Z_2$, $F/(E\cap F)\cong\Z_2$, the group
$\sg{E/(E\cap F), F/(E\cap F)}$ is isomorphic to $D_{2m}$ for some $m\geq 1$
(see e.g. \cite{CM}). Furthermore, since $$[\sg{E,F}:(E\cap F)]=|\sg{E,F}/(E\cap F)|$$ we have $[\sg{E,F}:(E\cap F)]=2m$.
Finally, it is clear that $$L(E\cap F,\langle E,F\rangle)\cong L\left(\langle E,F\rangle/(E\cap F)\right)$$ and therefore
$L(E\cap F,\sg{E,F})\cong L(D_{2m}).$
\qed

In the rest of this subsection it is assumed that $G\leq\sym(\Omega)$ is a transitive
permutation group containing $H$.

\bp \la{lp-} The lattice $L(G_\omega, G)$ is lower semimodular.
\ep
\pr
Assume the contrary and let
$E_1\in L(G_\omega,G)$ be a subgroup of $G$ such that \be \la{iiuu} E\cap F<E_1<E,\ee
where $E, F\in L(G_\omega, G),$ $E\neq F,$ are maximal in $\sg{E,F}$.
Notice that then
$$E_1\cap F = E\cap F.$$ If $E_1$ is permutable with $F$,
then \cmnt{$\sg{E_1,F} = E_1F$ and} by \eqref{form}
$$\left[\sg{E_1,F}:F\right]= \left[E_1:E_1\cap F\right]=\left[E_1:E\cap F\right]<
\left[E:E\cap F\right]\leq \left[\sg{E,F}:F\right].$$
Therefore, $\sg{E_1,F} < \sg{E,F}$. Since $F\leq \sg{E_1,F}$ and $F$ is maximal in
$\sg{E,F}$, this implies that $\sg{E_1,F} = F.$ Hence, $E_1\leq F$ and therefore
$E_1\leq E\cap F$ in contradiction with the assumption that $E\cap F < E_1.$

Suppose now that $F$ and $E_1$ are not permutable. Then Proposition \ref{properties} implies
that both $E_1$ and $F$ are not core-complementary. It follows now from
Propositions \ref{normality} and \ref{normal1} that there exist
$F',E_1'\in L_c(G_\omega,G)$ such that $[E_1:E'_1]=[F:F']=2$. Notice that each of the groups $F'$ and $E_1'$ is permutable with any $X\in L(G_\omega,G)$ by Proposition \ref{properties}. In particular,
$E_1'F\in L(G_\omega,G)$ and $EF'  \in L(G_\omega,G)$.

It follows from $$F\leq E_1'F \leq \sg{E,F}$$ that
either $E_1'F=\sg{E, F}$ or $E_1'F=F$. If $E_1'F=\sg{E, F}$, then the inclusions
$$
\sg{E,F}\supseteq EF\supseteq E_1F\supseteq E_1'F=\sg{E, F}
$$
imply that $E_1 F =\sg{E, F}\in L(G_\omega,G)$ in contradiction with the assumption that $E_1$ and $F$ are not permutable. So, assume that $E_1'F=F$. In this case
$E_1'\leq F$ and hence $E_1'\leq E\cap F.$
Together with $E\cap F< E_1$ and $[E_1:E_1']=2$ this implies that
\be \la{tre} E_1'=E\cap F = E_1\cap F.\ee In view of Proposition \ref{normal1} the last equality yields that
$E\cap F$ is $H$-intransitive and $E\cap F\leq F'$. Consequently, \be \la{posl} E\cap F = E\cap F'.\ee

It follows from $$E\leq EF' \leq \sg{E,F}$$ that either $EF'=\sg{E, F}$ or $EF'=E$.
If the equality $EF'=\sg{E, F}$ holds, then \eqref{form} and \eqref{posl} imply the inequality
$$
[F:E\cap F]\leq [\sg{E, F}:E]=[EF':E] = [F': E\cap F'] = [F':E\cap F] = \frac{1}{2} [F: E\cap F]$$ which is impossible. So, assume that $EF'=E$. In this case $F'\leq E$ and therefore $F'\leq E\cap F \leq F$.
Together with $[F:F']=2$ this implies that either $E\cap F = F$ or $E\cap F =F'$.
Furthermore, since in view of the maximality of $F$ and $E$ in $\sg{E, F}$
the equality $F\cap E = F$ is impossible, we may assume that $F'=E\cap F$. In this case
$[F:E\cap F]=2$. Together with \eqref{tre} and $[E_1:E_1'] = 2$ this implies that
\be \la{ppoo}
[F:E_1\cap F]=[E_1:E_1\cap F]=2.
\ee
It follows now from Proposition \ref{jk}
that the lattice $L(E_1\cap F,\sg{E_1, F})$ is isomorphic to the subgroup lattice of a
dihedral group $D_{2m}$, where $2m=[\sg{E_1, F}:E_1\cap F]$.
Furthermore, it follows from $F\leq \sg{E_1,F} \leq \sg{E,F}$ that either $\sg{E_1,F} = F$ or
$\sg{E_1,F} = \sg{E,F}$.
The first case is impossible since $E\cap F < E_1 < E$. Therefore
$\sg{E_1,F} = \sg{E,F}$ and hence \be \la{olo} L(E\cap F,\sg{E,F}) = L(E_1\cap F,\sg{E_1,F})\cong L(D_{2m}).\ee

Since maximal subgroups of $D_{2m}$ have prime index, it follows from \eqref{olo} that
the number $p:=[\sg{E,F}:F]$ is prime and hence
$$[\sg{E,F}:E\cap F]= [\sg{E,F}:F][F:E\cap F]=2p.$$
On the other hand, by \eqref{ppoo} $$[\sg{E,F}:E\cap F]=[\sg{E,F}:E][E:E_1][E_1:E\cap F] = 2[\sg{E,F}:E][E:E_1].$$
Therefore, $[\sg{E,F}:E][E:E_1]=p.$ Since this equality implies that at least one of the  numbers
$[\sg{E,F}:E],$ $[E:E_1]$ is equal to one, we conclude that there exists no $E_1\in L(G_\omega,G)$ satisfying
\eqref{iiuu} and therefore the lattice $L(G_\omega, G)$ is lower semimodular.
\qed

\bp\la{fin} Let $E,F\in L(G_\omega,G),$ $E\neq F$. Suppose that $E\cap F$ is maximal in $E$, $F$.
Then either $E$ and $F$ are permutable and $E$, $F$
are maximal in $\sg{E, F}$, or $E\cap F\normaleq \sg{E,F}$ and $\sg{E,F}/(E\cap F)\cong D_{2m}$ for some $m\geq 1$. Furthermore, $L(E\cap F,\sg{E,F})\cong L(D_{2m}).$

\ep
\pr
If $E$ and $F$ are permutable, then $E$ and $F$ are maximal in $\sg{E,F}=EF$
by Proposition \ref{lp}. So, suppose that $E$ and $F$ are not permutable and consider
the core-complementary subgroups $E'<E,$ $F'<F$ from Proposition \ref{normal1}.

It follows from $$E\cap F\leq E'(E\cap F)\leq E$$ that either $E'(E\cap F)=E$ or $E'(E\cap F)=E\cap F$.
In the first case we obtain $$EF = E'(E\cap F)F=E'F\in L(G_\omega,G)$$ that contradicts to the assumption that $E$ and $F$
are not permutable. Therefore $E' (E\cap F) = E\cap F$, or, equivalently, $E'\leq E\cap F$.
Since $[E:E']=2$, this implies that $E'=E\cap F$. Analogously, $F'=E\cap F$.
Thus $$[E:E\cap F] =[F:E\cap F] =2.$$ Now
Proposition~\ref{jk} yields the result.
\qed

\bc \la{gru}
Let $E,F\in L(G_\omega, G)$ be maximal in $\sg{E,F}$. Then
$E\cap F$ is maximal in $E$ and $F$ and either $EF=FE$, or $E\cap F\normaleq\sg{E,F}$
and $\sg{E,F}/(E\cap F)\cong D_{2m}$ for a prime $m$.
\ec
\pr  By Proposition~\ref{lp-} the group $E\cap F$ is maximal in $F$ and $E$. If $E$ and $F$ are
not permutable, then Proposition~\ref{fin} implies that $E\cap F\normaleq\sg{E,F}$
and $\sg{E,F}/(E\cap F)\cong D_{2m}$ for some $m\geq 1$.
Furthermore, since $F$ is maximal in $\sg{E,F}$ the group
$F/(E\cap F)\cong \Z_2$ is maximal in the group $\sg{F,E}/(E\cap F)\cong D_{2m}$ and therefore $m$ is prime.
\qed

We can summarize Propositions \ref{lp-} and \ref{fin} as follows.

\bt \la{bi} Let $G$ be a transitive permutation group containing a cyclic subgroup with two orbits.
Then the lattice
$L(G_\omega,G)$ is lower semi-modular.
Furthermore, $L(G_\omega,G)$ is modular unless there exists an interval
of $L(G_\omega,G)$ which
is isomorphic to the subgroup lattice of a dihedral group.
\et

\pr By Proposition~\ref{lp-} the lattice $L(G_\omega,G)$ is lower semi-modular.
If it is not modular, then existence of an interval isomorphic to $L(D_{2m})$ follows
from Proposition~\ref{fin}.

\bc \la{bi2} Let $F$ be a rational function such that its monodromy group contains a cyclic subgroup with at most
two orbits. Then any two maximal decompositions of $F$ are weakly equivalent. Furthermore, if
$$F=F_1\circ F_2\circ \dots \circ F_k \ \ \ \ {\it and} \ \ \ \
F=R_1\circ R_{2}\circ \dots \circ R_k$$ are two decompositions of $F$ then the set of degrees of the functions
$F_i,$ $1\leq i \leq k,$ coincides with the set of degrees of the functions
$G_i,$ $1\leq i \leq k.$
\ec

\pr
The first part of corollary follows from
Theorem \ref{bi} and Corollary \ref{df}. Furthermore, it follows
from the first part that in order to prove the second part
it is enough to establish that if $A,B$ are subgroups of $G$ such that
$A\cap B$ is maximal in $A,B$, and $A,B$ are maximal in $\sg{A,B}$
then the sets $\{[\sg{A,B}:B],[B:A\cap B]\}$ and
$\{[\sg{A,B}:A],[A:A\cap B]\}$ coincide.
If $A$ and $B$ are permutable, then this is a corollary of
formula \eqref{form}.
On the other hand, if $A$ and $B$ are not permutable then the property
needed easily follows from Corollary \ref{gru}.
\qed

\vskip 0.2cm
\noindent{\bf Remark.} The proof of Theorem \ref{bi} given above is a simplified
version of the proof given in the earlier preprint of the authors \cite{mp1}. Notice that Corollary \ref{bi2} and a weaker version of Corollary \ref{gru} were also independently proved in the
preprint \cite{kz} appeared shortly after \cite{mp1}.

\subsection{Non-permutable subgroups of $L(G_1,G)$ and algebraic curves having a factor
of genus zero with at most two points at infinity}
The following result is the algebraic counterpart of Proposition 2 in \cite{f2} (see also Theorem 8.1 of \cite{bilu} and Theorem 3.5 of \cite{pak}).

\bp\la{fri} Let $G$ be a group and $A,B$ be non-permutable subgroups of $G$.
Then there exist non-permutable subgroups $\hat A,\hat B$ of $G$ such that $A\leq \hat A,$ $B\leq \hat B,$
and $\core_G\hat A=\core_G\hat B.$
\ep
\pr
For $C\leq G$ denote by $d(C)$ a maximal number such that there
exists a maximal chain of subgroups
$$C=C_0< C_1<\dots < C_{d(C)}=G.$$
We use the induction on the number $d=d(A)+d(B)$. In order to lighten notation set $N=\core_GA,$ $M=\core_GB.$

First of all notice that the subgroups $AM$ and $BN$ are not permutable
since $$(AM)(BN)=AB, \ \ \ (BN)(AM)=BA.$$
In particular,
$AM\neq G$ and $BN\neq G$. So, if $d=2$ (that is if both $A$ and $B$ are maximal in $G$),
then $AM = A$, $BN=B$, and hence
$M\leq A$ and $N\leq B$. Since $M\normaleq G$ and $N\normaleq G$, this imples that
$M\leq N $ and $N\leq M$,
and hence $M=N.$
Therefore, if $d=2$ we can set $\hat{A}:=A,$ $\hat{B}:=B$.

Assume now that $d > 2$. If $d(AM) < d(A)$ or
$d(BN) < d(B)$, then the proposition follows from
the induction assumption. On the other hand, if $d(AM)=d(A)$ and $d(BN)=d(B)$, then
as above $AM = A$, $BN=B$, and $M=N.$ Therefore, we can set $\hat{A}:=A,\hat{B}:=B$.
\qed

Proposition \ref{fri} together with previous results allows us to describe
non-permu\-table subgroups of $L(G_\omega,G)$.

\bt \la{gru0} Let $G$ be a transitive permutation group containing a cyclic subgroup with two orbits and
$E,F\in L(G_\omega, G)$ be non-permutable subgroups of $G$ such that $\sg{E,F}=G$.
Then there exists $N\normaleq G$ such that
$E\cap F\leq N$
and $G/N \cong D_{2m}$ for some $m\geq 1$.
\et
\pr
By Proposition \ref{fri} there exist non-permutable subgroups
$\hat E,\hat F$ of $G$ such that $E\leq \hat E,$ $F\leq \hat F,$ and
$\core_G\hat E=\core_G \hat F.$  Furthermore, Proposition \ref{properties} implies
that both $\hat E$ and $\hat F$ are not core-complementary. Therefore,
by Propositions \ref{normal} and \ref{normal1}
\be \la{derf} [\hat E:\hat{E}']= 2, \ \ [\hat F:\hat{F}']= 2,\ee where
$\hat{E}' = (\core_G \hat E) G_\omega$ and $\hat{F}' = (\core_G \hat F) G_\omega$.

Since $\core_G\hat E=\core_G \hat F$, we obtain $\hat{E}'=\hat{F}'\leq \hat{E}\cap\hat{F}$.
On the other hand, the inequality $\hat E\hat F\neq \hat F \hat E$ implies that $\hat{E}\cap \hat{F}$ is a proper subgroup
of both $\hat{E}$ and $\hat{F}$. It follows now from \eqref{derf} that $\hat{E}'=\hat{F}' = \hat{E}\cap\hat{F}$ and
$[\hat{E}:\hat{F}\cap\hat{F}] = [\hat{F}:\hat{F}\cap\hat{F}] = 2$. Therefore, the theorem follows from
Proposition \ref{jk} taking into account that $E\cap F\leq \hat E\cap \hat F$.
\qed

Theorem \ref{gru0} has an interesting connection with the problem of description
of algebraic curves
\be \la{an2} A(x)-B(y)=0\ee
having a factor of genus zero with at most two points at infinity.
This problem is closely related to the number theory and in this context
was studied in the papers \cite{f1}, \cite{bilu}. In particular, in \cite{bilu}
a complete classification of such curves (defined over any field $k$ of characteristic zero) was obtained.
Another proof of this classification (over $\C$) was given in the paper \cite{pak} in the context of description of
double decompositions $$L=A\circ B=C\circ D$$
of rational functions $L$, with at most two poles, into compositions of rational functions.
The last problem turns out to be more general than the previous one since if curve
\eqref{an2} has an irreducible factor of genus zero with two points at infinity
then this factor may be parametrized by some Laurent polynomials and therefore there exist Laurent polynomials $L,L_1,L_2$ such that
the equality \be \la{dde} L=A\circ L_1=B\circ L_2\ee holds.

The both proofs of the classification of curves \eqref{an2} having a factor
of genus $0$ with at most two points at infinity split into two parts: the first one is the analysis of
the condition that, under the assumption that \eqref{an2} is irreducible, the genus
of \eqref{an2} is zero, and the second one is the reduction of the general case to the case when \eqref{an2} is irreducible. The first part essentially consists of
 a straightforward although highly laborious analysis of the formula which
calculates the genus of \eqref{an2} via the branching data of $A$ and $B$,
while the second part requires some more sophisticated considerations.

Denote by $G$ the monodromy group of $L$ and let $G_A$, $G_B$ be subgroups
of $L(G_\omega,G)$ corresponding to decompositions \eqref{dde}.
Then the condition that
\eqref{an2} is reducible is equivalent to the condition that
$G_AG_B\neq G.$
Therefore, Theorem \ref{gru0} can be viewed as
an algebraic counterpart of the portion of the discussed classification related to the reducible case,
and implies easily the corresponding result
(cf. \cite{bilu}, Theorem 9.3 and \cite{pak}, Theorem 7.3).

\bp \la{red}
Suppose that
curve \eqref{an2} is reducible and has a factor of genus zero with at most two points at infinity. Then
there exist polynomials $R,$ $\tilde A,$ $\tilde B,$ $\mu,$ where $\deg \mu=1,$
such that
\be \la{ggvv} A=R \circ \tilde A,\ \ \  B =R\circ \tilde B\ee
and either the curve $\tilde A(x)- \tilde B(y)=0$ is irreducible,
or \be \la{vvgg} \tilde A=- T_{lr}\circ \mu, \ \ \ \tilde B=T_{ls}\circ \mu,\ee where $T_{lr}, T_{ls}$
are the corresponding Chebyshev polynomials with $r,s\geq 1$, $l>2,$ and $\GCD(r,s)=1.$
\ep
\pr Without loss of generality we may assume that there exists no polynomial
$R,$ $\deg R>1,$ such that \eqref{ggvv} holds
for some polynomials $\tilde A,\tilde B$, or equivalently that $\sg{G_A,G_B}=G.$
If curve \eqref{an2} is irreducible, then there is nothing to prove so assume that \eqref{an2} is reducible.
In this case $L_1,$ $L_2$ are not polynomials since otherwise Corollary \ref{eng}
and the assumption about solutions of \eqref{ggvv} imply the equality $\gcd(\deg A,\deg B)=1$ which in its turn implies easily the
irreducibility of curve \eqref{an2}.
Therefore, the cyclic subgroup $H$ of $G$ generated by the permutation corresponding to a loop around infinity has two orbits.

It follows now from Theorem \ref{gru0} that there exists
$N\normaleq G$ such that $N\in L(G_\omega,G)$
and $G/N \cong D_{2m}$ for some $m\geq 1$.
Furthermore, since $N\normaleq G$ the action of $G$ on cosets of $N$ is regular.
Therefore, $$G//N\cong G/N \cong D_{2m}$$ and hence there exists
a decomposition $L=U\circ V$ of $L$ such that the monodromy group of $U$ is a regular
covering of the sphere with the dihedral monodromy group. By the well known
classification of regular coverings of the sphere
which goes back to Klein (see \cite{klein} and the Appendix below) this implies that
$$U=\mu_1\circ \frac{1}{2}\left(z^{m}+\frac{1}{z^{m}}\right)\circ \mu_2,$$
where $\mu_1,\mu_2$ are automorphisms of the sphere.

Clearly, without loss
of generality we may assume that $\mu_1=z$.
Furthermore, since $L$ has poles only at the points $0$ and $\infty$
it follows from $L=U\circ V$
that $\mu_2 \circ V=z^{\pm n}\circ (cz)$ for some $n\geq 1$ and $c\in \C$. Therefore, \be \la{erty} L=\frac{1}{2}\left(z^{mn}+\frac{1}{z^{mn}}\right)\circ (cz)\ee
and $G=D_{2mn}.$
Now the proposition follows easily from the description of possible double decompositions
of function \eqref{erty}. \qed

\section{Appendix}
In this appendix we describe the structure of maximal decompositions of
rational functions which are
regular coverings of the sphere that is of the functions for which $G_{\omega}=e.$
These functions, appearing in a variety of different contexts from differential equations to Galois theory,
were first described by Klein in \cite{klein}.
For such a function $f$ its monodromy group $G$ is isomorphic to its
automorphism group and therefore is isomorphic to a finite subgroup of $\Aut\C\P^1.$ Any such
a subgroup is isomorphic to one of the groups $C_n,$ $D_{2n},$ $A_4,$ $S_4,$ $A_5$ and the corresponding function $f$ is defined by its group
up to a composition $\mu_1\circ f\circ \mu_2,$ where $\mu_1,$ $\mu_2 \in \Aut\C\P^1$. 

The Klein functions provide the simplest examples of rational functions for which
the first Ritt theorem fails to be true. Indeed, if $f$ is a Klein function then
its maximal decompositions
correspond to maximal chains of subgroups
of its monodromy group $G.$ Therefore, in order to find counterexamples to the first Ritt
theorem it is enough to find non $r$-equivalent maximal chains of subgroups of $G$. For the groups $C_n$ and $D_n$ such chains do not exist while
for the groups $A_4,$ $S_4,$ $A_5$ they do. For example, it is easy to see that
\be \la{chch} e < C_2 <   V_4 < A_4,\ \ \ e < C_3 < A_4,\ee where $C_2$ (resp. $C_3$) is a cyclic group of order 2 (resp. 3) and $V_4$ is the Klein
four group, are maximal chains of different length in $A_4$ and therefore
for the corresponding Klein function the first Ritt theorem
fails to be true.
The fact that the first Ritt theorem is not true for arbitrary rational functions
was observed already by
Ritt itself in \cite{r1}. Although Ritt did not give any indications about the nature of such examples (see
the discussion in \cite{g1}, \cite{g2}, \cite{be1}), the fact that the Klein functions
corresponding to $A_4,$ $S_4,$ $A_5$
were mentioned by him in a close context in the paper \cite{r2} suggests that he meant exactly these functions.

Below we give a detailed analysis of decompositions
of the Klein functions.
We show that
for a function $f$ corresponding to $A_4$ or $S_4$ the number of
weak equivalence classes of its maximal decompositions equals two and that
two non-equivalent maximal decompositions of $f$ are
weakly equivalent if and only if they have the same length. On the other hand, we show that
the function corresponding to $A_5$ has six weak equivalence classes of maximal decompositions
five of which have the same length.
Besides, we give several related explicit examples of non weakly equivalent maximal decompositions.
In particular, we
give an example of a rational function with three poles
for which the first Ritt theorem fails to be true.

\subsection{Decompositions of $f_{C_n}$ and $f_{D_{2n}}$}\label{meloch}
For the
cyclic and dihedral groups the representatives of the corresponding classes of
Klein functions are $$f_{C_n}=z^n, \ \ \ \ \ \ \ \ f_{D_{2n}}=\frac{1}{2}\left(z^n+\frac{1}{z^n}\right)$$ and
by Corollary \ref{bi2} all maximal decompositions of these functions are weakly equivalent.
Observe that any decomposition of $f_{C_n}$ into a composition of two functions is equivalent to the decomposition $$z^{n/d}\circ z^d,$$ where $d\vert n,$ while any decomposition of $f_{D_{2n}}$ is equivalent either
to the decomposition
$$\frac{1}{2}\left(z^n+\frac{1}{z^n}\right)=\frac{1}{2}\left(z^{n/d}+\frac{1}{z^{n/d}}\right)\circ z^d,$$
where $d\vert n,$
or to the decomposition $$
\frac{1}{2}\left(z^{n}+\frac{1}{z^{n}}\right)=\mu^{n/d} T_{n/d}\circ \frac{1}{2}\left(\mu z^{d}+\frac{1}{\mu z^{d}}\right),
$$ where $d\vert n$ and $\mu^{2n/d}=1.$

\subsection{Decompositions of $f_{A_4}$} The subgroup lattice of the group $A_4$ can be described as follows.
$A_4$ has tree subgroups $C_2^1,C_2^2,C_2^3$ of order 2 which are conjugated between themselves and are contained in a unique subgroup of order 4 which
is the Klein four-group $V_4=\{e,(12)(34),(13)(24),(14)(23)\}$. Besides, $A_4$ has
four conjugated subgroups $C_3^1,C_3^2,C_3^3,C_3^4$ of order 3 which are maximal in $A_4$.
This implies that
$f_{A_4}$ has 7 non-equivalent decompositions corresponding to the chains
\be \la{ch1} e < C_2^1 <  V_4 < A_4, \ \ \ \ e < C_2^2 <   V_4 < A_4,
\ \ \ \ e < C_2^3 <   V_4 < A_4,\ee and
\be \la{ch2} e < C_3^1 < A_4, \ \ \ \ e < C_3^2 < A_4, \ \ \ \ e < C_3^3 < A_4, \ \ \ \ e < C_3^4 < A_4.\ee

Clearly, all decompositions from the first group are
$r$-equivalent. The same is true for decompositions from the second group.
On the other hand, compositions from the first and the second groups obviously are non-equivalent
since they have different lengths.

\subsection{Decompositions of $f_{S_4}$} Similarly to the case of the group $A_4$
two maximal chains in $S_4$ are $r$-equivalent if and only if they have the same length. However, since $S_4$ has already 28 proper subgroups, in order to prove this statement we will use an argument distinct from the examination of all maximal chains.

First of all, notice that any maximal subgroup of
$S_4$ either is $A_4$, or is conjugate to
$$D_8=\{e, (12),(34),(12)(34),(13)(24),(14)(23),(1324), (1432)\},$$ or is conjugate to
$S_3$. 
Besides,
it is easy to see that any maximal chain of subgroups of $A_4$ has length $3$ or $4$.
We show now that any two maximal chains
$$\f F:\ 1 < F_1 < F_2 < S_4 \ \ \ \ {\rm and} \ \ \ \ \f E:\  1 < E_1 < E_2 < S_4$$
of length $3$ are $r$-equivalent.
If $E_2=F_2$, then the statement is clear so we
may assume that $E_2\neq F_2.$ This implies
in particular that $E_2\cap F_2$ is a proper subgroup of the groups $E_2$ and $F_2$.
Observe that $E_2\cap F_2$ is non-trivial
since otherwise we would have $|S_4|\geq |E_2||F_2|\geq  36>\vert S_4\vert $.
In order to prove that the chains $\f F$ and $\f E$ 
are $r$-equivalent it is enough to show that the chains
$$\tilde{\f F}:\ 1 < F_2\cap E_2 < F_2 < S_4 \ \ \ \ {\rm and} \ \ \ \ \tilde{\f E}:\  1 < F_2\cap E_2 < E_2 < S_4$$
are maximal since then
$$
\f F \sim \tilde{\f F} \sim \tilde{\f E}
\sim \f E.
$$
First, notice that $E_2, F_2\not\cong D_8$, since maximal chains
in $D_8$ have length 3. Therefore,
at least
one of the groups $E_2$, $F_2$, say $F_2$, is isomorphic
$S_3$ and hence the chain $\tilde{\f F}$ is maximal since $\vert S_3\vert =6$.
If $E_2\cong S_3$, then the chain $\tilde{\f E}$ is maximal as well.
On the other hand, if $E_2 = A_4$, then
$\vert F_2\cap E_2\vert =\vert S_3\cap A_4\vert =3$ implying that
the chain $1 < F_2\cap E_2 < E_2$
is one of the chains \eqref{ch2} and, therefore, is maximal.

Similarly, any two chains
$$\f F:\ 1 < F_1 < F_2 < F_3 <S_4 \ \ \ \ {\rm and} \ \ \ \ \f E:\  1 < E_1 < E_2 < E_3 < S_4$$
of length $4$ are $r$-equivalent. Indeed, if $E_3=F_3$ then either
$E_3 = F_3\cong D_8$ or $E_3 = F_3 = A_4$ and the statement is true since maximal chains of
equal length in the groups $D_8$ and $A_4$ are $r$-equivalent.
Therefore, we may assume that $F_3=A_4$,
$E_3 = D_8.$ Setting now
$$V_4=\{e,(12)(34),(13)(24),(14)(23)\},\ \ \ C_2= \{e,(12)(34)\}$$
and observing that
$E_3\cap F_3 = V_4$, we see that
the chains
$$\tilde{\f F}:\ 1< C_2 < V_4 < A_4 < S_4, \ \ \ \ {\rm and} \ \ \ \ \tilde{\f E}:\ 1 < C_2 < V_4 < D_8<S_4$$
are maximal.
Since any two chains of equal length inside $D_8$ and $A_4$ are equivalent, this implies that
$$\f F \sim \tilde{\f F} \sim \tilde{\f E}
\sim \f E.
$$
\vskip 0.2cm

\subsection{Decompositions of $f_{A_5}$} It is easy to see that any maximal subgroup of $A_5$ is conjugated
either to $A_4$, or to $D_{10}$, or to $S_3$ and that any maximal chain of subgroups in $f_{A_5}$ has length
3 or 4.
In contrast to the groups $A_4$, $S_4$ in the group $A_5$ we face a new phenomenon:
although any two maximal chains of length $3$
in $A_5$ are $r$-equivalent there exist non $r$-equivalent decompositions of length
$4$.

First prove that any two maximal chains $$\f F:\ 1 < F_1 < F_2 < A_5 \ \ \ \
{\rm and} \ \ \ \ \f E:\  1 < E_1 < E_2 < A_5$$ of length $3$ in $A_5$ are $r$-equivalent.
If $E_2=F_2$, then the statement is clear so we
may suppose that $E_2\neq F_2.$

Assume first that $E_2\cong D_{10}$ and $F_2\cong S_3$. Since $A_5$
is not a product of $D_{10}$ and $S_3$, the intersection $E_2\cap F_2$ is non-trivial. Therefore
the chains
$$\tilde{\f F}:\ 1 < F_2\cap E_2 < F_2 < A_5
\ \ \ \ {\rm and} \ \ \ \ \tilde{\f E}:\  1 < F_2\cap E_2 < E_2 < A_5$$ are maximal,
implying
$$\f F \sim \tilde{\f F} \sim \tilde{\f E}\sim \f E. $$
By transitivity of $\sim$ this yields that
any two maximal chains of length $3$ such that $E_2\cong S_3,$ $F_2\cong S_3$ or $E_2\cong D_{10},$ $F_2\cong D_{10}$ also are $r$-equivalent.

Let now $$\f B:\  1< B_1< B_2<A_5$$
be a maximal chain such that
$B_2\cong A_4$. Then \eqref{ch2} implies that $|B_1|=3$. One can check that  
the normalizer $C$ of any group of order 3 in $A_5$ is isomorphic to $S_3$.
Therefore, $\f B$ is
equivalent to a maximal chain
$$1 < B_1 < C < A_5$$ with $C\cong S_3$. It follows now from the transitivity of $\sim$
that all the chains of length $3$ are $r$-equivalent.

Let us show now that two maximal chains of length $4$
$${\f B}:=1 < B_1 < B_2 < B_3 < A_5 \ \ \ \ {\rm and} \ \ \ \ {\f C}:=1 < C_1 < C_2 < C_3 < A_5$$
in $A_5$ are equivalent if and only if their maximal subgroups
coincide. Clearly, we have $B_3\cong C_3 \cong A_4$. If $B_3 = C_3$, then $\f B \sim \f C$
since any two chains of length 4 in $A_4$ are $r$-equivalent.

Assume now that $B_3\neq C_3$.
If the chains ${\f B}$ and ${\f C}$ are equivalent, then in the sequence of maximal chains which connects them
there should be two chains of the form $$1 < P_1 < P_2 < P_3 < A_5, \ \ \ \ \ 1 < P_1 < P_2 < Q_3 < A_5,$$
where $P_3\neq Q_3$. The maximality condition implies that $P_3\cap Q_3=P_2$. Furthermore,
$P_2 \cong V_4$ by \eqref{ch1}. On the other hand,
$A_4$ contains a unique Sylow $2$-subgroup of order $4$ which is normal
in $A_4$. Therefore, $P_2\trianglelefteq P_3$, $P_2\trianglelefteq Q_3$ and hence
$P_2\trianglelefteq \sg{P_3,Q_3}=A_5.$ Since this contradicts to the simplicity of $A_5$, we conclude that
${\f B}$ and ${\f C}$ are not $r$-equivalent. \qed

\subsection{Explicit formulas} Although all the information about maximal decompositions
of Klein functions can be obtained from the analysis given above,
the actual finding of the corresponding decompositions requires some non trivial calculations.
In particular, the corresponding
maximal decompositions which do not satisfy the first Ritt theorem were found explicitly
only for the simplest chains \eqref{chch} (see \cite{be2}, \cite{gs}). It turns out that a convenient tool for such calculations is the Grothendieck theory of ``Dessins d'enfants" which
provides an identification of $f_{A_4},$ $f_{S_4},$ and $f_{A_5}$ with the Belyi functions of
the tetrahedron, cube, and octahedron respectively.
Below we give several explicit examples of non equivalent maximal decompositions
obtained by this method, referring the reader
interested in details of calculations to the paper \cite{pazv}.

First, a calculation shows that the Belyi functions for the tetrahedron
can be written in the form
\be \la{tet} f_{A_4}=-\frac{1}{64}\frac{z^3(z^3-8)^3}{(z^3+1)^3}\ee
and any maximal decomposition of $f_{A_4}$ is weakly equivalent either to
$$f_{A_4}=-\frac{1}{64}\frac{z(z-8)^3}{(z+1)^3}\circ z^3$$ or to
the decomposition
\be \la{tet1}
f_{A_4}=-\frac{1}{64}z^3\circ \frac{z^2-4}{z-1}\circ \frac{z^2+2}{z+1}.
\ee

Furthermore, one can show that the inclusion $A_4\subset S_4$ implies that
\be \la{s5} f_{S_4}= -\frac {4x}{{x}^{2}+1-2\,x}\circ f_{A_4}=
\frac{256{z}^{3} \left( {z}^{6}-7z^3-8 \right) ^{3}}
{\left( z^6 +20z^3-8\right) ^{4}}\ee
and therefore the decompositions of $f_{S_4}$ corresponding to the chains
$$1 < C_3 < A_4 < S_4, \ \ \ \ \ 1 < C_2 < V_4 < A_4 < S_4$$
are
$$f_{S_4}= \left(-\frac {4x}{{x}^{2}+1-2\,x}\right)\circ\left( -\frac{1}{64}\frac{z(z-8)^3}{(z+1)^3}\right)\circ z^3,$$
and $$f_{S_4}= \left(-\frac {4x}{{x}^{2}+1-2\,x}\right)\circ \left(-\frac{1}{64}z^3\right)\circ \left(\frac{z^2-4}{z-1}\right)\circ \left(\frac{z^2+2}{z+1}\right).$$

On the other hand, one can show that
for example the maximal decompositions of $f_{S_4}$ (written in a bit different normalization) corresponding to the chains
\be \la{sx} 1 < C_2 < C_4 < D_8 < S_4,  \ \ \ \ \ 1 < C_2 < S_3 < S_4\ee
are:
$$-\frac{1}{432}\frac{(16x^8-56x^4+1)^3}{x^4(4x^4+1)^4}=
\left(\frac{1}{54}\frac{(z+7)^3}{(z-1)^2}\right)\circ \left(\frac{1}{2}\left(z+\frac{1}{z}\right)\right)\circ (-z^2)\circ 2z^2$$ and
$$-\frac{1}{432}\frac{(16x^8-56x^4+1)^3}{x^4(4x^4+1)^4}=\left(-\frac{256}{27}z^3(z-1)\right)\circ \left(\frac{1}{4}\frac{(z-1)^3}{z^2+1}+1\right)\circ
\left(z-\frac{1}{2z}\right).$$

Finally, identifying the chains of subgroups
\be \la{konez} C_2<S_3<S_4, \ \ \ C_2<V_4<D_8<S_4\ee
with maximal decompositions of
the function \be \la{fufuf} -\frac{1}{27}\frac{(z^4+2z^2-3)^3}{(z^2+1)^4}\ee which is
a left compositional factor of $f_{S_4}$,
one can show that to \eqref{konez} correspond the maximal decompositions:
$$-\frac{1}{27}\frac{(z^4+2z^2-3)^3}{(z^2+1)^4}=
\left(\frac{1}{54}\frac{(7-z)^3}{(z+1)^2}\right)\circ \left(2z^2+4z+1\right)\circ z^2
$$
and
$$-\frac{1}{27}\frac{(z^4+2z^2-3)^3}{(z^2+1)^4}=\left(-\frac{256}{27}z^3(z-1)\right)\circ \left(\frac{1}{4}\frac{(z-1)^3}{z^2+1}+1\right).$$
Notice that since function \eqref{fufuf} has three poles this example shows that with no additional
assumptions the first Ritt theorem can not be extended to rational functions the monodromy of which
contains a cyclic subgroup with more than two orbits.

\vskip 0.2cm 
\noindent{\bf Remark.} It is interesting to understand how wide is the class of
rational function for which the
first Ritt theorem holds. The following observation can be useful
for obtaining some experimental results in this direction.

Suppose that $F$ is a rational function providing a counterexample to the first Ritt theorem
and let $N_F$ be the Galois closure of $F$.
Then $N_F$ is a regular covering of the sphere and $F$ is a left compositional factor of $N_F$.
In particular, the functions entering
into decompositions of $F$ are left
factors of $N_F.$
Therefore, all possible counterexamples to the first Ritt theorem can be obtained by the analysis of  regular coverings only.

This observation suggests to analyse possible counterexamples to the first Ritt theorem
relatively to the genus of their Galois closure.
For example, the rational functions for which this genus equals zero are
exactly the Klein functions
and their compositional left factors considered above. The next case which would be interesting
to investigate is the one
corresponding to rational factors of regular covering of genus 1.

Notice also that the analysis of decompositions of $f_{A_5}$ suggests that an other source of possible counterexamples
to the first Ritt theorem is the functions which admit non-equivalent decompositions of the
form $A\circ B=A\circ D$. Examples of such functions can be found in the papers \cite{r2}, \cite{az}, \cite{zv}, \cite{pazv}.

\bibliographystyle{amsplain}

\begin{thebibliography}{10}

\bibitem {ai} M. Aigner,
\textit{Combinatorial theory},
Grundlehren der mathematischen Wissenschaften 234, Springer-Verlag, (1979).

\bibitem {az} R. Avanzi, U. Zannier,
\textit{The equation $f(X) = f(Y)$ in rational functions $X = X(t)$, $Y = Y(t)$},
Compos. Math. 139, No. 3, 263-295 (2003)


\bibitem {be1} W. Bergweiler, \textit{An example concerning factorization of rational functions}, Exposition. Math.  11  (1993),  no. 3, 281--283.

\bibitem {be2} W. Bergweiler, \textit{Erratum to the paper ``An example concerning factorization of rational functions"}, http://analysis.math.uni-kiel.de/bergweiler/schrift.html{\#}preprints.



\bibitem {bilu} Y. Bilu, R. Tichy,
\textit{The Diophantine equation $f(x) = g(y)$},
Acta Arith. 95, No.3, 261-288 (2000).


\bibitem {cg} J. M. Couveignes, L. Granboulan, \textit{Dessins from a geometric point of view},
in The Grothendieck theory of dessins d'enfants, 79-113. Cambridge University Press, 1994.

\bibitem{CM} Coxeter H.S.M., Moser W.O.J.,
\textit{Generators and Relations for Discrete Groups,} 4th ed., Springer Verlag, New York, 1980.

\bibitem{dm} J.D. Dixon and B. Mortimer, \textit{Permutation groups,} Springer, 1996.

\bibitem {en} H. Engstrom, \textit{Polynomial substitutions,} Amer. J. Math. 63, 249-255 (1941).

\bibitem {f1} M. Fried,
\textit{On a theorem of Ritt and related diophantine problems},
J. Reine Angew. Math. 264, 40-55 (1973).

\bibitem {f2} M. Fried,
\textit{Fields of definition of function fields and a problem in the reducibility of polynomials in two variables}, Ill. J. Math. 17, 128-146 (1973).


\bibitem {f3} M. Fried, \textit{
Arithmetical properties of function fields. II. The generalized Schur problem,} Acta Arith. 25, 225-258 (1974)


\bibitem {g1} F. Gross, \textit{On factorization theory of meromorphic functions},  Comment. Math. Univ. St. Paul.  24, no. 1, 47--60 (1975/76).

\bibitem {g2} F. Gross, \textit{Factorization of meromorphic functions and some open problems}, in Complex analysis, Lecture Notes in Math., Vol. 599, pp. 50-67, Springer, Berlin, 1977.

\bibitem {gs} J. Gutierrez, D. Sevilla,
\textit{Building counterexamples to generalizations for rational functions of Ritt's decomposition theorem},
J. Algebra 303, No. 2, 655-667 (2006).

\bibitem{huppert} B. Huppert, \textit{ Endliche Gruppen, I.} Springer-Verlag, 1967.

\bibitem {klein} F. Klein,
\textit{Lectures on the icosahedron and the solution of equations of the fifth degree},
New York: Dover Publications, (1956).


\bibitem{kz} G. Kuperberg, M. Zieve,\textit{
Analogues of the Jordan-H{\" o}lder theorem for transitive
$G$-sets}, preprint, arXiv:0712.4142v1.


\bibitem {kur} A. Kurosch, \textit{
The theory of groups. Vol. I,} New York: Chelsea Publishing Company, (1955).

\bibitem {zv} N. Magot, A. Zvonkin, \textit{
Belyi functions for Archimedean solids},
Discrete Math. 217, No.1-3, 249-271 (2000).

\bibitem {mu} P. M\"uller, \textit{Primitive monodromy groups of polynomials},
in ``Recent developments in the inverse Galois problem.'', Providence, RI: American Mathematical Society. Contemp. Math. 186, 385-401 (1995).

\bibitem {mz} P. M{\" u}ller, M.Zieve, \textit{On Ritt's polynomial decomposition theorem}, preprint, arXiv:0807.3578v1.

\bibitem {mp1} M. Muzychuk, F. Pakovich, \textit{On maximal decompositions of rational functions}, preprint, arXiv:0712.3869v1

\bibitem{p1} F. Pakovich, \textit{On polynomials sharing preimages of compact sets, and related questions}, Geom. Funct. Anal, 18, No. 1, 163-183 (2008).

\bibitem{p2} F. Pakovich, \textit{On the functional equation $F(A(z))=G(B(z))$, where $A,B$ are polynomial and $F,G$ are continuous functions}, Math. Proc. Cambridge Philos. Soc. 143, No.2, 469-472,
(2007).

\bibitem{pppp} F. Pakovich, \textit{The algebraic curve $P(x)-Q(y)=0$ and
functional equations,} Complex Var. Elliptic Equ., to appear, arXiv:0804.0736v2.

\bibitem {ar} F. Pakovich, {\it On analogues of Ritt theorems for rational functions with at most two poles}, Russ. Math. Surv. 63, No. 2, 181-182 (2008).

\bibitem{pak} F. Pakovich, \textit{Prime and composite Laurent polynomials}, preprint arXiv:0710.3860v5.

\bibitem{pazv} F. Pakovich, A. Zvonkin, \textit{Dessins d'enfants and functional equations},
in preparation.

\bibitem {r1} J Ritt,
\textit{Prime and composite polynomials},
American M. S. Trans. 23, 51-66 (1922).

\bibitem {r2} J Ritt,
\textit{Equivalent rational substitutions},
American M. S. Trans. 26, 221-229 (1924).

\bibitem {r3} J Ritt,
\textit{Permutable rational functions,}
American M. S. Trans. 25, 399-448 (1923).

\bibitem {sch} A. Schinzel, \textit{Polynomials with special regard to
reducibility}, Encyclopedia of Mathematics and Its Applications
77, Cambridge University Press, 2000.

\bibitem {tor} P. Tortrat, \textit{Sur la composition des polyn\^omes},
Colloq. Math. 55, No.2, 329-353 (1988).



\bibitem{wi} H. Wielandt,
\textit{Finite Permutation Groups},
 Academic Press, 1964, Berlin.


\bibitem{zi} M. Zieve, \textit{Decompositions of Laurent polynomials}, preprint, arXiv:0710.1902v1.


\end{thebibliography}

\end{document}